

\ifx\begin\undefined\else\global\advance\srcdepth by
1\expandafter \fi

\def\begin{}
\newcount\srcdepth
\srcdepth=1

\outer\def\bye{\global\advance\srcdepth by -1
  \ifnum\srcdepth=0
    \def\endcmd{\vfill\eject\nopagenumbers\par\vfill\supereject\end}
  \else\def\endcmd{}\fi
  \endcmd
}


\magnification=\magstephalf
\baselineskip=13pt
\hsize = 5.5truein
\hoffset = 0.5truein
\vsize = 8.5truein
\voffset = 0.2truein
\emergencystretch = 0.05\hsize

\newif\ifblackboardbold

\blackboardboldtrue


\font\titlefont=cmbx12 scaled\magstephalf
\font\sectionfont=cmbx12


\newfam\bboldfam
\ifblackboardbold
\font\tenbbold=msbm10
\font\sevenbbold=msbm7
\font\fivebbold=msbm5
\textfont\bboldfam=\tenbbold
\scriptfont\bboldfam=\sevenbbold
\scriptscriptfont\bboldfam=\fivebbold
\def\bbold{\fam\bboldfam\tenbbold}
\else
\def\bbold{\bf}
\fi


\font\Arm=cmr8
\font\Ai=cmmi8
\font\Asy=cmsy8
\font\Abf=cmbx8
\font\Brm=cmr6
\font\Bi=cmmi6
\font\Bsy=cmsy6
\font\Bbf=cmbx6
\font\Crm=cmr5
\font\Ci=cmmi5
\font\Csy=cmsy5
\font\Cbf=cmbx5

\ifblackboardbold
\font\Abbold=msbm10 at 8pt
\font\Bbbold=msbm7 at 6pt
\font\Cbbold=msbm5
\fi

\def\smallmath{%
\textfont0=\Arm \scriptfont0=\Brm \scriptscriptfont0=\Crm
\textfont1=\Ai \scriptfont1=\Bi \scriptscriptfont1=\Ci
\textfont2=\Asy \scriptfont2=\Bsy \scriptscriptfont2=\Csy
\textfont\bffam=\Abf \scriptfont\bffam=\Bbf \scriptscriptfont\bffam=\Cbf
\def\rm{\fam0\Arm}\def\mit{\fam1}\def\oldstyle{\fam1\Ai}%
\def\bf{\fam\bffam\Abf}%
\ifblackboardbold
\textfont\bboldfam=\Abbold
\scriptfont\bboldfam=\Bbbold
\scriptscriptfont\bboldfam=\Cbbold
\def\bbold{\fam\bboldfam\Abbold}%
\fi
}








\newlinechar=`@
\def\forwardmsg#1#2#3{\immediate\write16{@*!*!*!* forward reference should
be: @\noexpand\forward{#1}{#2}{#3}@}}
\def\nodefmsg#1{\immediate\write16{@*!*!*!* #1 is an undefined reference@}}

\def\forwardsub#1#2{\def\newref{{#2}{#1}}}

\def\forward#1#2#3{%
\expandafter\expandafter\expandafter\forwardsub\expandafter{#3}{#2}
\expandafter\ifx\csname#1\endcsname\relax\else%
\expandafter\ifx\csname#1\endcsname\newref\else%
\forwardmsg{#1}{#2}{#3}\fi\fi%
\expandafter\let\csname#1\endcsname\newref}

\def\firstarg#1{\expandafter\argone #1}\def\argone#1#2{#1}
\def\secondarg#1{\expandafter\argtwo #1}\def\argtwo#1#2{#2}

\def\ref#1{\expandafter\ifx\csname#1\endcsname\relax
  {\nodefmsg{#1}\bf`#1'}\else
  \expandafter\firstarg\csname#1\endcsname
  ~\expandafter\secondarg\csname#1\endcsname\fi}

\def\refs#1{\expandafter\ifx\csname#1\endcsname\relax
  {\nodefmsg{#1}\bf`#1'}\else
  \expandafter\firstarg\csname #1\endcsname
  s~\expandafter\secondarg\csname#1\endcsname\fi}

\def\refn#1{\expandafter\ifx\csname#1\endcsname\relax
  {\nodefmsg{#1}\bf`#1'}\else
  \expandafter\secondarg\csname #1\endcsname\fi}



\def\widow#1{\vskip 0pt plus#1\vsize\goodbreak\vskip 0pt plus-#1\vsize}



\def\marginlabel#1{}

\def\showlabelsabove{
\font\labelfont=cmss10 at 6pt
\def\marginlabel##1{\rlap{\smash{\raise 10pt\hbox{\labelfont##1}}}}
}

\newcount\seccount
\newcount\proccount
\seccount=0
\proccount=0

\def\stdskip{\vskip 9pt plus3pt minus 3pt}
\def\stdbreak{\par\removelastskip\penalty-100\stdskip}

\def\proof{\stdbreak\noindent{\bf Proof. }}

\def\qed{\vrule height 1.2ex width .9ex depth .1ex}

\def\Box{
  \ifmmode\eqno\qed
  \else\ifvmode\removelastskip\line{\hfil\qed}
  \else\unskip\quad\hskip-\hsize
    \hbox{}\hskip\hsize minus 1em\qed\par
  \fi\stdbreak\fi}

\def\references{
  \removelastskip
  \widow{.05}
  \vskip 24pt plus 6pt minus 6 pt
  \leftline{\sectionfont References}
  \nobreak\stdskip\noindent}

\def\ifempty#1#2\endB{\ifx#1\endA}
\def\makeref#1#2#3{\ifempty#1\endA\endB\else\forward{#1}{#2}{#3}\fi}

\outer\def\section#1 #2\par{
  \removelastskip
  \global\advance\seccount by 1
  \global\proccount=0\relax
                \edef\numtoks{\number\seccount}
  \makeref{#1}{Section}{\numtoks}
  \widow{.05}
  \vskip 24pt plus 6pt minus 6 pt
  \message{#2}
  \leftline{\marginlabel{#1}\sectionfont\numtoks\quad #2}
  \nobreak\stdskip}

\def\proclamation#1#2{
  \outer\expandafter\def\csname#1\endcsname##1 ##2\par{
  \stdbreak
  \advance\proccount by 1
  \edef\numtoks{\number\seccount.\number\proccount}
  \makeref{##1}{#2}{\numtoks}
  \noindent{\marginlabel{##1}\bf #2 \numtoks\enspace}
  {\sl##2\par}
  \stdbreak}}

\def\othernumbered#1#2{
  \outer\expandafter\def\csname#1\endcsname##1{
  \stdbreak
  \advance\proccount by 1
  \edef\numtoks{\number\seccount.\number\proccount}
  \makeref{##1}{#2}{\numtoks}
  \noindent{\marginlabel{##1}\bf #2 \numtoks\enspace}}}

\proclamation{definition}{Definition}
\proclamation{lemma}{Lemma}
\proclamation{proposition}{Proposition}
\proclamation{theorem}{Theorem}
\proclamation{corollary}{Corollary}
\proclamation{conjecture}{Conjecture}

\othernumbered{example}{Example}
\othernumbered{remark}{Remark}
\othernumbered{construction}{Construction}





\input epsf

\newcount\figcount
\figcount=0
\newbox\drawing
\newcount\drawbp
\newdimen\drawx
\newdimen\drawy
\newdimen\ngap
\newdimen\sgap
\newdimen\wgap
\newdimen\egap

\def\drawbox#1#2#3{\vbox{
  \setbox\drawing=\vbox{\offinterlineskip\epsfbox{#2.eps}\kern 0pt}
  \drawbp=\epsfurx
  \advance\drawbp by-\epsfllx\relax
  \multiply\drawbp by #1
  \divide\drawbp by 100
  \drawx=\drawbp truebp
  \ifdim\drawx>\hsize\drawx=\hsize\fi
  \epsfxsize=\drawx
  \setbox\drawing=\vbox{\offinterlineskip\epsfbox{#2.eps}\kern 0pt}
  \drawx=\wd\drawing
  \drawy=\ht\drawing
  \ngap=0pt \sgap=0pt \wgap=0pt \egap=0pt 
  \setbox0=\vbox{\offinterlineskip
    \box\drawing \ifgridlines\drawgrid\drawx\drawy\fi #3}
  \kern\ngap\hbox{\kern\wgap\box0\kern\egap}\kern\sgap}}

\def\draw#1#2#3{
  \setbox\drawing=\drawbox{#1}{#2}{#3}
  \advance\figcount by 1
  \goodbreak
  \midinsert
  \centerline{\ifgridlines\boxgrid\drawing\fi\box\drawing}
  \smallskip
  \vbox{\offinterlineskip
    \centerline{Figure~\number\figcount}
    \smash{\marginlabel{#2}}}
  \endinsert}

\def\nextfigtoks{%
  \advance\figcount by 1%
  \edef\numtoks{\number\figcount}%
  \advance\figcount by -1}

\newif\ifgridlines
\newbox\figtbox
\newbox\figgbox
\newdimen\figtx
\newdimen\figty

\newdimen\bwd
\bwd=2sp 

\def\hline#1{\vbox{\smash{\hbox to #1{\leaders\hrule height \bwd\hfil}}}}

\def\vline#1{\hbox to 0pt{%
  \hss\vbox to #1{\leaders\vrule width \bwd\vfil}\hss}}

\def\clap#1{\hbox to 0pt{\hss#1\hss}}
\def\vclap#1{\vbox to 0pt{\offinterlineskip\vss#1\vss}}

\def\hstutter#1#2{\hbox{%
  \setbox0=\hbox{#1}%
  \hbox to #2\wd0{\leaders\box0\hfil}}}

\def\vstutter#1#2{\vbox{
  \setbox0=\vbox{\offinterlineskip #1}
  \dp0=0pt
  \vbox to #2\ht0{\leaders\box0\vfil}}}

\def\crosshairs#1#2{
  \dimen1=.002\drawx
  \dimen2=.002\drawy
  \ifdim\dimen1<\dimen2\dimen3\dimen1\else\dimen3\dimen2\fi
  \setbox1=\vclap{\vline{2\dimen3}}
  \setbox2=\clap{\hline{2\dimen3}}
  \setbox3=\hstutter{\kern\dimen1\box1}{4}
  \setbox4=\vstutter{\kern\dimen2\box2}{4}
  \setbox1=\vclap{\vline{4\dimen3}}
  \setbox2=\clap{\hline{4\dimen3}}
  \setbox5=\clap{\copy1\hstutter{\box3\kern\dimen1\box1}{6}}
  \setbox6=\vclap{\copy2\vstutter{\box4\kern\dimen2\box2}{6}}
  \setbox1=\vbox{\offinterlineskip\box5\box6}
  \smash{\vbox to #2{\hbox to #1{\hss\box1}\vss}}}

\def\boxgrid#1{\rlap{\vbox{\offinterlineskip
  \setbox0=\hline{\wd#1}
  \setbox1=\vline{\ht#1}
  \smash{\vbox to \ht#1{\offinterlineskip\copy0\vfil\box0}}
  \smash{\vbox{\hbox to \wd#1{\copy1\hfil\box1}}}}}}

\def\drawgrid#1#2{\vbox{\offinterlineskip
  \dimen0=\drawx
  \dimen1=\drawy
  \divide\dimen0 by 10
  \divide\dimen1 by 10
  \setbox0=\hline\drawx
  \setbox1=\vline\drawy
  \smash{\vbox{\offinterlineskip
    \copy0\vstutter{\kern\dimen1\box0}{10}}}
  \smash{\hbox{\copy1\hstutter{\kern\dimen0\box1}{10}}}}}

\def\figtext#1#2#3#4#5{
  \setbox\figtbox=\hbox{#5}
  \dp\figtbox=0pt
  \figtx=-#3\wd\figtbox \figty=-#4\ht\figtbox
  \advance\figtx by #1\drawx \advance\figty by #2\drawy
  \dimen0=\figtx \advance\dimen0 by\wd\figtbox \advance\dimen0 by-\drawx
  \ifdim\dimen0>\egap\global\egap=\dimen0\fi
  \dimen0=\figty \advance\dimen0 by\ht\figtbox \advance\dimen0 by-\drawy
  \ifdim\dimen0>\ngap\global\ngap=\dimen0\fi
  \dimen0=-\figtx
  \ifdim\dimen0>\wgap\global\wgap=\dimen0\fi
  \dimen0=-\figty
  \ifdim\dimen0>\sgap\global\sgap=\dimen0\fi
  \smash{\rlap{\vbox{\offinterlineskip
    \hbox{\hbox to \figtx{}\ifgridlines\boxgrid\figtbox\fi\box\figtbox}
    \vbox to \figty{}
    \ifgridlines\crosshairs{#1\drawx}{#2\drawy}\fi
    \kern 0pt}}}}


\def\hpad#1#2#3{\hbox{\kern #1\hbox{#3}\kern #2}}
\def\vpad#1#2#3{\setbox0=\hbox{#3}\dp0=0pt\vbox{\kern #1\box0\kern #2}}



\def\stack#1#2#3{\vbox{\offinterlineskip
  \setbox2=\hbox{#2}
  \setbox3=\hbox{#3}
  \dimen0=\ifdim\wd2>\wd3\wd2\else\wd3\fi
  \hbox to \dimen0{\hss\box2\hss}
  \kern #1
  \hbox to \dimen0{\hss\box3\hss}}}


\def\hexp#1{%
  \setbox0=\hbox{${}^{#1}$}%
  \hbox to .5\wd0{\box0\hss}}



\def\bmatrix#1#2{{\smallmath\left[\vcenter{\halign
  {&\kern#1\hfil$##\mathstrut$\kern#1\cr#2}}\right]}}

\def\rightarrowmat#1#2#3{
  \setbox1=\hbox{\kern#2$\bmatrix{#1}{#3}$\kern#2}
  \,\vbox{\offinterlineskip\hbox to\wd1{\hfil\copy1\hfil}
    \kern 3pt\hbox to\wd1{\rightarrowfill}}\,}

\def\leftarrowmat#1#2#3{
  \setbox1=\hbox{\kern#2$\bmatrix{#1}{#3}$\kern#2}
  \,\vbox{\offinterlineskip\hbox to\wd1{\hfil\copy1\hfil}
    \kern 3pt\hbox to\wd1{\leftarrowfill}}\,}

\def\rightarrowbox#1#2{
  \setbox1=\hbox{\kern#1\hbox{\smallmath #2}\kern#1}
  \,\vbox{\offinterlineskip\hbox to\wd1{\hfil\copy1\hfil}
    \kern 3pt\hbox to\wd1{\rightarrowfill}}\,}

\def\leftarrowbox#1#2{
  \setbox1=\hbox{\kern#1\hbox{\smallmath #2}\kern#1}
  \,\vbox{\offinterlineskip\hbox to\wd1{\hfil\copy1\hfil}
    \kern 3pt\hbox to\wd1{\leftarrowfill}}\,}






\def\bookletdims{
  \hsize=5.25truein
  \vsize=7truein
}

\def\legalbooklet#1{
  \input quire
  \bookletdims
  \htotal=7.0truein
  \vtotal=8.5truein
  \hoffset=\htotal
  \advance\hoffset by -\hsize
  \divide\hoffset by 2
  \voffset=\vtotal
  \advance\voffset by -\vsize
  \divide\voffset by 2
  \advance\voffset by -.0625truein
  \shhtotal=2\htotal
  \horigin=0.0truein
  \vorigin=0.0truein
  \shstaplewidth=0.01pt
  \shstaplelength=0.66truein
  \shthickness=0pt
  \shoutline=0pt
  \shcrop=0pt
  \shvoffset=-1.0truein
  \ifnum#1>0\quire{#1}\else\qtwopages\fi
}

\def\preview{
  \input quire
  \bookletdims
  \hoffset=0.1truein
  \vtotal=8.5truein
  \shhtotal=14truein
  \voffset=\vtotal
  \advance\voffset by -\vsize
  \divide\voffset by 2
  \advance\voffset by -.0625truein
  \htotal=2\hoffset
  \advance\htotal by \hsize
  \horigin=0.0truein
  \vorigin=0.0truein
  \shstaplewidth=0.5pt
  \shstaplelength=0.5\vtotal
  \shthickness=0pt
  \shoutline=0pt
  \shcrop=0pt
  \shvoffset=-1.0truein
  \qtwopages
}

\def\twoup{
  \input quire
  \hsize=4.79452truein 
  \vsize=7truein
  \vtotal=8.5truein
  \shhtotal=11truein
  \hoffset=-2\hsize
  \advance\hoffset by \shhtotal
  \divide\hoffset by 6
  \voffset=\vtotal
  \advance\voffset by -\vsize
  \divide\voffset by 2
  \advance\voffset by -12truept
  \htotal=2\hoffset
  \advance\htotal by \hsize
  \horigin=0.0truein
  \vorigin=0.0truein
  \shstaplewidth=0.01pt
  \shstaplelength=0pt
  \shthickness=0pt
  \shoutline=0pt
  \shcrop=0pt
  \shvoffset=-1.0truein
  \qtwopages
}


\newcount\countA
\newcount\countB
\newcount\countC

\def\monthname{\begingroup
  \ifcase\number\month
    \or January\or February\or March\or April\or May\or June\or
    July\or August\or September\or October\or November\or December\fi
\endgroup}

\def\dayname{\begingroup
  \countA=\number\day
  \countB=\number\year
  \advance\countA by 0 
  \advance\countA by \ifcase\month\or
    0\or 31\or 59\or 90\or 120\or 151\or
    181\or 212\or 243\or 273\or 304\or 334\fi
  \advance\countB by -1995
  \multiply\countB by 365
  \advance\countA by \countB
  \countB=\countA
  \divide\countB by 7
  \multiply\countB by 7
  \advance\countA by -\countB
  \advance\countA by 1
  \ifcase\countA\or Sunday\or Monday\or Tuesday\or Wednesday\or
    Thursday\or Friday\or Saturday\fi
\endgroup}

\def\timename{\begingroup
   \countA = \time
   \divide\countA by 60
   \countB = \countA
   \countC = \time
   \multiply\countA by 60
   \advance\countC by -\countA
   \ifnum\countC<10\toks1={0}\else\toks1={}\fi
   \ifnum\countB<12 \toks0={\sevenrm AM}
     \else\toks0={\sevenrm PM}\advance\countB by -12\fi
   \relax\ifnum\countB=0\countB=12\fi
   \hbox{\the\countB:\the\toks1 \the\countC \thinspace \the\toks0}
\endgroup}

\def\timestamp{\dayname, \the\day\ \monthname\ \the\year, \timename}


\def\enma#1{{\ifmmode#1\else$#1$\fi}}

\def\mathbb#1{{\bbold #1}}
\def\mathbf#1{{\bf #1}}

\def\NN{\enma{\mathbb{N}}}
\def\ZZ{\enma{\mathbb{Z}}}



\def\bb{\enma{\mathbf{b}}}
\def\cc{\enma{\mathbf{c}}}
\def\xx{\enma{\mathbf{x}}}
\def\mm{\enma{\mathbf{m}}}

\def\th{{^{\rm th}}}

\def\H{\widetilde{H}}
\def\mtext#1{\;\,\allowbreak\hbox{#1}\allowbreak\;\,}

\def\set#1{\enma{\{#1\}}}
\def\setdef#1#2{\enma{\{\;#1\;\,|\allowbreak
  \;\,#2\;\}}}

\def\idealdef#1#2{\enma{\langle\;#1\;\,
  |\;\,#2\;\rangle}}

\def\depth{\mathop{\rm depth}\nolimits}
\def\reg{\mathop{\rm reg}\nolimits}

\def\Ext{\mathop{\rm Ext}\nolimits}
\def\Tor{\mathop{\rm Tor}\nolimits}
\def\link{\mathop{\rm lk}\nolimits}
\def\star{\mathop{\rm star}\nolimits}
\def\core{\mathop{\rm core}\nolimits}
\def\abs#1{\enma{\left| #1 \right|}}

\def\GCD{\mathop{\rm GCD}\nolimits}

\def\In{\mathop{\rm In}\nolimits}
\def\Gin{\mathop{\rm Gin}\nolimits}

\def\red{\mathop{\rm red}\nolimits}

\def\Ker{\mathop{\rm Ker}\nolimits}
\def\Hom{\mathop{\rm Hom}\nolimits}
\def\Syz{\mathop{\rm Syz}\nolimits}

\def\onedot{\drawbox{50}{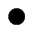}{}}
\def\twodots{\drawbox{50}{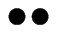}{}}
\def\threedots{\drawbox{50}{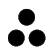}{}}
\def\fourdots{\drawbox{50}{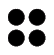}{}}
\def\fivedots{\drawbox{50}{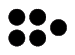}{}}

\def\onedota{\lower 1pt \drawbox{50}{onedot}{}}
\def\twodotsa{\lower 1pt \drawbox{50}{twodots}{}}
\def\threedotsa{\lower 1pt \drawbox{50}{threedots}{}}
\def\fourdotsa{\lower 1pt \drawbox{50}{fourdots}{}}
\def\fivedotsa{\lower 1pt \drawbox{50}{fivedots}{}}

\hbox{}
\bigskip
\centerline{\titlefont Extremal Betti Numbers
and 
Applications to Monomial Ideals} \bigskip
\centerline{Dave Bayer \quad Hara Charalambous \quad Sorin Popescu
\footnote{$^{*}$}{\rm The first and third authors are grateful to the NSF for 
support during the preparation of this work. The second author is grateful 
to Columbia University for its hospitality during the 
preparation of this manuscript.}}
\bigskip
\bigskip

Let $S=k[x_1,\dots,x_n]$ be the polynomial ring in $n$ 
variables over a field $k$, let $M$ be a graded $S$-module, and let 
$$F_\bullet: \quad 0  \longrightarrow F_r 
\longrightarrow \ldots \longrightarrow 
F_1 \longrightarrow  F_0 \longrightarrow M \longrightarrow 0$$ 
be a  minimal free resolution of $M$ over $S$. As usual, we define 
the  associated (graded) Betti numbers $\beta_{i,j}=\beta_{i,j}(M)$ 
by the formula
$$
F_i = \oplus_j S(-j)^{\beta_{i,j}}.
$$
Recall that the (Mumford-Castelnuovo) regularity of $M$
is the least integer $\rho$ such that for each $i$
all free generators of $F_i$ lie in degree $\le i+\rho$, that
is $\beta_{i,j}=0$, for $j>i+\rho$.
In terms of {\it Macaulay} [Mac] regularity is the number of rows in 
the diagram produced by the ``betti'' command.

A Betti number $\beta_{i,j}\ne 0$ will be called {\sl extremal}
if $\beta_{l,r}=0$ for all $l\ge i$, $r\ge j+1$ and $r-l\ge j-i$, 
that is if $\beta_{i,j}$ is the nonzero top left ``corner'' 
in a block of  zeroes in the {\it Macaulay} ``betti'' diagram. 
In other words, extremal Betti numbers account for
``notches'' in the shape of the minimal free resolution
and one of them computes the regularity. In this sense, 
extremal Betti numbers can be seen as a refinement of the
notion of Mumford-Castelnuovo regularity.

In the first part of this note we connect the extremal Betti numbers 
of an arbitrary submodule of a free $S$-module with those of its 
generic initial module. 
In the second part, which can be read independently of the first, we relate
extremal multigraded Betti numbers in the minimal resolution of a square 
free monomial ideal with those of the monomial ideal corresponding
to the Alexander dual simplicial complex.

Our techniques give also a simple geometric proof of a more precise version 
of a recent result of Terai [Te97] (see also [FT97] for a 
homological reformulation and related results), 
generalizing Eagon and Reiner's theorem 
[ER96] that a Stanley-Reisner ring is Cohen-Macaulay if and only 
if the homogeneous ideal corresponding to the Alexander dual simplicial 
complex has a  linear resolution.

We are grateful to David Eisenbud for 
useful discussions.

\section{gins} GINs and extremal Betti numbers

A theorem of Bayer and Stillman ([BaSt87], [Ei95])
asserts that if $M$ is a graded submodule of a free 
$S=k[x_1,\dots,x_n]$-module
$F$, and one considers the degree reverse lexicographic
monomial order, then after a generic change of coordinates, 
the modules $F/M$ and $F/\In(M)$ have the same regularity and the same
depth (in this situation the module $\In(M)$ is known as $\Gin(M)$).

We generalize  this result to show that  
corners in the minimal resolution of $F/M$ correspond to
corners in the minimal resolution of  $F/\Gin(M)$  
and that moreover the extremal Betti numbers of $F/M$
and of $F/\Gin(M)$ match. The proof, inspired by 
the approach in [Ei95],  shows 
that each extremal Betti number of
$F/M$ or respectively $F/\Gin(M)$ is computed by the 
unique extremal Betti number of a finite length 
submodule.

We use the same notation as in the introduction. 
Let $M$ be a graded $S$-module, and let 
$$F_\bullet: \quad 0  \longrightarrow F_r 
\longrightarrow \ldots \longrightarrow 
F_1 \longrightarrow  F_0 \longrightarrow M \longrightarrow 0$$ 
be a  minimal free resolution of $M$. As usual, we define  
$\Syz_l(M):=\Ker( F_l\longrightarrow F_{l-1})$ to be the 
$l^{\rm{th}}$ syzygy module of $M$.
 
We say that $M$ is {\it $(m,l)$-regular} iff $\Syz_l(M)$ is
$(m+l)$-regular (in the classical sense); that is to say that 
all generators of $F_j$ for $l\le j\le r$
have degrees $\le j+m$. 

We also define the
{\it $l$-regularity} of $M$, denoted in the sequel as $l$-$\reg(M)$, 
to be the regularity of the module
$\Syz_l(M)$; it is the least integer
$m$ such that $M$ is $(m,l)$-regular. 

It is easy to see that 
$\reg(M)=0$-$\reg(M)$ and $l$-$\reg(M)\le (l-1)$-$\reg(M)$.
Strict inequality occurs only at extremal Betti numbers,
which thus pinpoint ``jumps'' in the regularity of the 
successive syzygy modules. In this case, 
if $m=l$-$\reg(M)$, we say that 
$(l,m)$ is a {\it corner} of $M$ and
that $\beta_{l,m+l}(M)$ is an {\it extremal Betti number} 
of $M$.

\proposition{propbA} $M$ is $(m,l)$-regular iff
$$\Ext^j(M,S)_k=0\quad \hbox{for all $j\ge l$ and all $k\le -m-j-1$}.$$
If moreover $(l,m)$ is a corner of $M$, then 
$\beta_{l,m+l}(M)$ is equal to the 
number of minimal generators of  $\Ext^l(M,S)$ in degree $(-m-l)$.

\proof The first part follows from [Ei95, Proposition 20.16] 
since  $M$ is $(m,l)$-regular iff $\Syz_l(M)$ is
$m$-regular. For the second part notice
that by degree considerations the nonzero generators of 
$F_l$ in degree $m+l$
correspond to nontrivial cycles of $\Ext^l(M,S)_{-m-l}$.\Box

Finite length modules have exactly one
extremal Betti number:

\theorem{thmbA} If $M$ is a finite length module,
and $\beta_{l,m+l}$ is an extremal Betti number of $M$, then
$l=n$ and $\beta_{n,m+n}$
is the last  nonzero value in
the Hilbert function of $M$.

\proof Since $M$ has finite length it follows that 
$\Ext^j(M,S)=0$ for all $j<n$. On the other hand
$\Ext^n(M,S)_{-n-t}= \Ext^n(M,S(-n))_{-t}
\cong \Hom_k (M,k)_{-t}= \Hom_k(M_t,k)$,  
from which the claim follows easily.
\Box 

\corollary{corbA} Let $F$ be a graded free 
$S$-module with basis, and let $M$ a graded submodule
of $F$ such that $F/M$ has finite length.
Then the extremal Betti number of $F/M$ is equal to the extremal
Betti number of $F/\Gin(M)$ (with respect to the graded reverse
lexicographic order). 

\proof $F/M$ and $F/\Gin(M)$ have the same Hilbert function.
\Box

\proposition{propbB} If $0\longrightarrow A\longrightarrow B
 \longrightarrow C \longrightarrow 0$ is a short exact sequence of graded 
finitely generated $S$-modules, then
\item{a)} ${l}$-$\reg(A)\le \max (l$-$\reg(B),(l+1)$-$\reg(C)+1)$.
\item{b)} $l$-$\reg(B)\le \max (l$-$\reg(A),l$-$\reg(C))$.
\item{c)} $l$-$\reg(C)\le \max ((l-1)$-$\reg(A)-1,l$-$\reg(B))$.
\item{d)} If $A$ has finite length, 
then  $l$-$\reg(B)= \max (l$-$\reg(A), l$-$\reg(C))$.

\proof The proof follows by examining the appropriate 
$l^{\rm{th}}$  graded pieces
of the long exact sequence in $\Ext(\cdot,S)$. 
See the analogue statement
for regularity in [Ei95, Corollary 20.19].
\Box

\lemma{lembA} If $F$ is a finitely generated graded free $S$-module,
$M$ a graded submodule of $F$, and $x$ a linear form of $S$  
such that the module $(M:x)/M$ 
has finite length, then 
$$l\hbox{\rm{-}$\reg$}(F/M)=\max (l\hbox{\rm{-}$\reg$}((M:x)/M),
l\hbox{\rm{-}$\reg$}(F/(M:x)) ).$$

\proof The claim follows from the short exact sequence
$$0\longrightarrow (M:x)/M \longrightarrow F/M
 \longrightarrow F/(M:x) \longrightarrow 0$$ and
\ref{propbB}.\Box

We may now state and prove  the analogue of Bayer and Stillman's 
([BaSt87], [Ei95]) result on regularity:

\theorem{thmbB} Let $F$ be a finitely generated 
graded free $S$-module with basis, let $M$ be a graded submodule
of $F$, and let $\beta_{i,j}$ denote the $i$th graded Betti number of $F/M$,
and $\beta^{gin}_{i,j}$ the $i$th graded Betti number of $F/\Gin(M)$. Then
$$l\hbox{\rm{-}$\reg$}(F/M)=l\hbox{\rm{-}$\reg$}(F/\Gin(M)).$$
If moreover $(l,m)$ is a corner of $F/M$, then
$$\beta_{l,m+l}=\beta^{gin}_{l,m+l}.$$

\proof We can assume that $\In(M)=\Gin(M)$. 
If $x_n$ is a nonzero divisor
of $F/M$ the claims follow by induction on the number of variables:
the Betti numbers  of $F/M$ over $S$ are equal to the Betti numbers of
$F/(x_nF,M)$ over $S/x_n$ and the initial module  of $M$ over $S$ is the
same as the initial module of $M/x_nM$ over $S/x_n$.
Therefore we will assume in the sequel 
that $x_n$ is a zero divisor of $F/M$.

We prove the first part of the theorem.
Since $(M:x_n)/M$
is a finite length module,
$n$-$\reg((M:x_n)/M)=\reg((M:x_n)/M)=\reg((\In(M):x_n)/\In(M))=
n$-$\reg((\In(M):x_n)/\In(M))$
and the first part of  the theorem follows by \ref{lembA} and 
induction on the sum of degrees of the elements in a reduced
Gr\"obner basis of $M$. (Recall that with $F$ and $M$ as above, 
using the reverse lexicographic order, if 
$\{g_1,\dots, g_t\}$ is a (reduced) Gr\"obner basis for $M$ and
$g_i':=g_i/\GCD(x_n,g_i)$, then $\{g_1',\dots, g_t'\}$ is a 
(reduced) Gr\"obner basis for the module $(M:x_n)$.)

Assume first that $(m,n)$ is a corner of $F/M$, so in particular
$\Ext^n(F/M,S)\neq 0$. Let $N=H^0_\mm(F/M)$ be the set of all elements
in $F/M$ that are annihilated by some power of the
ideal $\mm\subset S$ generated by the variables, and let $L:=(F/M)/N$.
From the short exact sequence
$$0\longrightarrow H^0_\mm(F/M) \longrightarrow F/M
 \longrightarrow L \longrightarrow 0,$$ we conclude that
$\Ext^n(F/M,S)\cong \Ext^n(H^0_\mm(F/M),S)$,
since $L$ has no torsion and thus
$\Ext^n(L,S)=0$.
By the first part and \ref{thmbA}
the last nonzero value, say $w$, of the
Hilbert function of $(M:x_n)/M$ (or $(\In(M):x_n)/\In(M)$)
occurs in degree $m$.
But this is also the last nonzero value of the
Hilbert function of $H^0_\mm(F/M)$. By \ref{corbA}
it follows that $\beta_{n,(m+n)}=w=\beta^{gin}_{n,(m+n)}$.

Finally consider a corner say $(m,l)$, with $l<n$, in the resolution of $F/M$.
From the short exact sequence in the proof of \ref{lembA}, it follows that
$\Ext^l(F/M,S)\cong \Ext^l(F/(M:x_n),S)$ so we are  done again
by induction on the sum of degrees of the elements in a reduced Gr\"obner
basis of $M$.\Box

Given an $S$-module $P$ define $\red(P)$ to be $P/(H^0_\mm(P)+xP)$,
where $x$ is a generic linear form. Let  $P_0=P$ and 
define $P_{i+1}=\red (P_i)$, for all $i\ge 1$.

\corollary{gen} Let $L$ be a module over
a polynomial ring $S$ with free presentation $L=F/M$,
and let $N=F/\Gin(M)$.  Then for all $i\geq 0$:
\item{a)} The Hilbert functions of $H^0_\mm(L_i)$ and $H^0_\mm(N_i)$
coincide,
\item{b)} The depths of $L_i$ and $N_i$ coincide.
\item{c)} The extremal Betti numbers of $L$ correspond to
jumps in the highest socle degrees of the $L_i$s.

\proof If $\depth (L)\geq 1$, then taking generic initial modules 
commutes with factoring out a generic linear form, 
up to semicontinuous numerical data such as depth and Hilbert function. 
The same thing is true if we factor out an element of highest degree
of the socle of $L$ since it corresponds to a corner by \ref{thmbB}. 
Induction now proves {\sl a)} and {\sl b)}.\Box

\section{alexander} Alexander duality and square-free monomial ideals

The minimal free resolution of a multigraded ideal in 
$S=k[x_1,\ldots,x_n]$, the polynomial ring in $n$  variables 
over a field $k$, is obviously multigraded, and so it is natural
to introduce and study in this context a multigraded  analogue for 
``extremal Betti numbers''.

We use the same notation as above. 
Let $S=k[x_1,\ldots,x_n]$ be the polynomial ring 
let $[n] = \set{1,\ldots,n}$, and let $\Delta$ denote the set of all subsets of
$[n]$. Given a simplicial complex $X\subseteq \Delta$, define the {\it
Stanley-Reisner\/} ideal $I_X\subseteq S$ to be the ideal generated by the
monomials corresponding to the nonfaces of $X$:
$$I_X \quad = \quad \idealdef{\xx^F}{F\not\in X}.$$
$I_X$ is a square-free monomial ideal, and every square-free monomial ideal 
arises in this way. 

Define the {\it Alexander dual\/} simplicial complex
$X^\vee\subseteq \Delta$ to be the complex obtained by 
successively complementing the faces of
$X$ and $X$ itself, in either order. In other words, define    
 $$ X^\vee \quad = \quad \setdef{F}{F^c \not\in X} 
\quad = \quad \Delta\setminus \setdef{F}{F^c \in X}$$ 
where $F^c$ denotes the complement $[n]\setminus F$.
Defining also the {\it Alexander dual\/} ideal $I_{X^\vee}$, note the
following pattern: $$\matrix{X & \longleftrightarrow & I_X \cr\cr
 \big\updownarrow & & \big\updownarrow \cr\cr
 I_{X^\vee} &  \longleftrightarrow & X^\vee \cr}$$
The sets of faces which define $X$, $X^\vee$, $I_X$ and $I_{X^\vee}$  are related
horizontally by complementing with respect to $\Delta$, and vertically by
complementing with respect to $[n]$.

The following is a simplicial version of Alexander duality:

\theorem{thmA}
 Let $X\subset \Delta$ be a simplicial complex. For any abelian group $G$, there are
isomorphisms
 $$\H_i(X;G) \;\cong\; \H^{n-i-3}(X^\vee;G) \; \mtext{and} \; \H^i(X;G)
\;\cong\; \H_{n-i-3}(X^\vee;G)$$ where $\H$ denotes reduced simplicial
(co)homology.

\proof
 First, suppose that $X$ is a nonempty, proper subcomplex of the sphere
$S^{n-2} = \Delta \setminus [n]$. Working with geometric
realizations, Alexander duality asserts (compare [Mun84, Theorem 71.1]) that
 $$\H_i(X;G) \;\cong\; \H^{n-i-3}(S^{n-2}\setminus X;G) \; \mtext{and} \;
\H^i(X;G) \;\cong\; \H_{n-i-3}(S^{n-2}\setminus X;G).$$
 The claim follows because $X^\vee$ is homotopy-equivalent to $S^{n-2}\setminus
X$: Let $X'$ denote the first barycentric subdivision of $X$. Complementing the faces
of $X^\vee$ embeds  $(X^\vee)'$ as a simplicial
subcomplex of $(S^{n-2}\setminus X)'$. The straight-line homotopy defined by
collapsing each face of $(S^{n-2}\setminus X)'$ 
onto its vertices not belonging to $X'$
is a strong deformation retract of $(S^{n-2}\setminus X)'$ onto $(X^\vee)'$.

The remaining cases $X=\emptyset$, $\set{\emptyset}$, $\Delta \setminus
[n]$, and $\Delta$ are easily checked by hand.\Box

\ref{thmA} is also easily proved directly, modulo a subtle sign change. Define a
pairing on faces $F$, $G\in\Delta$ by
 $$ \langle F,G\rangle \quad = \quad \cases{
 (-1)^{\lfloor{\abs F \over 2}\rfloor} \;\sigma(F,G), & if $G=F^c$ \cr
 0, & otherwise \cr}$$
 where $\sigma(F,G)$ is the sign of the permutation that  sorts  the concatenated
sequence $F$, $G$ into order. This pairing allows us to reinterpret any $i$-chain as
an $(n-i-2)$-cochain, identifying relative homology with relative cohomology. 
We compute 
 $$\H_i(X;G) \;\cong\; \H_{i+1}(\Delta,X;G) \;\cong\; \H^{n-i-3}(X^\vee,\emptyset;G)
\;\cong\; \H^{n-i-3}(X^\vee;G);$$ the second isomorphism is similar. 
See [Bay96] for
details. This formulation can also be understood as the self-duality of the Koszul
complex; see [BH93, 1.6.10].

\medskip

Given an arbitrary monomial ideal $I\subseteq S$, let
 $$L_\bullet: \quad 0  \longrightarrow L_m \longrightarrow \ldots \longrightarrow 
L_1 \longrightarrow  L_0 \longrightarrow I \longrightarrow 0$$
 be a minimal free resolution of $I$; we have $m \le n-1$. 
The multigraded {\it Betti numbers\/}
of $I$ are the ranks $\beta_{i,\bb} = \dim (L_i)_\bb$ of the 
$\bb\th$ graded summands
$(L_i)_\bb$ of $L_i$.

For each degree $\bb\in\NN^n$, define the following subcomplex of $\Delta$:
$$K_\bb(I) \quad = \quad \setdef{F\in\Delta}{\xx^{\bb-F} \in I}.$$
Here, we identify each face $F\in\Delta$ with its characteristic vector
$F\in\set{0,1}^n$. The following is a characterization of the Betti numbers of $I$ in
terms of $K_\bb(I):$

\theorem{thmB}
 The Betti numbers of a monomial ideal $I\subseteq S$ are given by
 $$\beta_{i,\bb} \quad = \quad \dim \H_{i-1}(K_\bb(I); k).$$

\proof
 The groups $\Tor_i(I,k)$ can be computed either by tensoring a resolution of $I$
by $k$, or by tensoring a resolution of $k$ by $I$. Using the minimal resolution 
$L_\bullet$ of $I$, one sees that $\beta_{i,\bb} = \dim \Tor_i(I,k)_\bb$.
Using the  Koszul complex $K_\bullet$ of $k$, $\Tor_i(I,k)$
is also the $i$th homology of the complex
 $$I\otimes K_\bullet: \quad 0  \longrightarrow I \otimes \wedge^n V
\longrightarrow \ldots \longrightarrow  I \otimes \wedge^1 V \longrightarrow 
I \otimes \wedge^0 V \longrightarrow 0,$$
 where $V$ is the subspace of degree one forms of $S$.
Now, $(I \otimes \wedge^i V)_\bb$ has a basis consisting of all expressions of
the form
 $$\xx^\bb/x_{j_1}\cdots x_{j_i} \; \otimes \; x_{j_1} \wedge \ldots \wedge x_{j_i}$$
 where $\xx^\bb/x_{j_1}\cdots x_{j_i} \in I$. These expressions correspond $1:1$ to
the $(i-1)$-faces $F=\set{j_1, \ldots, j_i}$ of $K_\bb(I)$. Thus, one recognizes $(I
\otimes K_\bullet)_\bb$ as the augmented oriented chain complex used to compute
$\H_{i-1}(K_\bb(I); k)$.
 \Box

A striking reformulation of \ref{thmB} for square-free monomial ideals is due to
Hochster [Ho77], based on ideas of Reisner [Rei76]. For each
$\bb\in\set{0,1}^n$, let $X_\bb$ denote the full subcomplex of $X$ on the vertices in
the support of \bb.

\theorem{thmC}
 Let $I_X\subseteq S$ be the square-free monomial 
ideal determined by the simplicial
complex $X\subseteq \Delta$.  We have $\beta_{i,\bb} = 0$ unless
$\bb\in\set{0,1}^n$, in which case
 $$\beta_{i,\bb} \quad = \quad \dim \H_{\abs\bb-i-2}(X_\bb; k).$$

\proof
 If $b_j>1$ for some $j$ then $K_\bb(I)$ is a cone over the vertex 
$j$, so $\beta_{i,\bb}= 0$ by \ref{thmB}. 
Otherwise, $X_\bb$ is the dual of $K_\bb(I)$ with respect to the
support of \bb:  $F\in K_\bb(I) \Leftrightarrow \xx^{\bb-F}\in I
\Leftrightarrow \bb-F \not\in X$. By \ref{thmA},
 $$\H_{i-1}(K_\bb(I); k) \quad \cong \quad \H^{\abs\bb-i-2}(X_\bb; k).$$
Homology and cohomology groups with coefficients in $k$ are
(non-canonically) isomorphic, so the result follows by \ref{thmB}.
 \Box

This is essentially Hochster's original argument; he implicitly proves Alexander
duality in order to interpret $\Tor_i(I,k)_\bb$ as computing the homology of $X_b$.

Recall that the {\it link\/} of a face $F\in X$ is the set
 $$
\link(F,X) \quad = \quad \setdef{G}{F\cup G \in X \mtext{and} F \cap G =
\emptyset}.
$$
 Together with the restrictions $X_\bb$, the links 
$\link(F,X)$ are the other key
ingredient in the study of square-free monomial ideals, 
dating to [Rei76]. They too
have a duality interpretation, first made explicit in [ER96]:
 The Betti numbers $\beta^\vee_{i,\bb}$ of $I_{X^\vee}$ can be 
computed using links in $X$.

For each $\bb\in\set{0,1}^n$, let $\bb^c$ denote the complement 
$(1,\ldots,1)-\bb$.

\theorem{thmD}
  Let $I_{X^\vee}\subseteq S$ be the square-free monomial ideal determined 
by the dual $X^\vee$ of the simplicial complex $X\subseteq \Delta$. We have
$\beta^\vee_{i,\bb} = 0$ unless $\bb\in\set{0,1}^n$ and 
$\bb^c \in X$, in which case
 $$\beta^\vee_{i,\bb} \quad = \quad \dim \H_{i-1}(\link(\bb^c,X); k).$$

\proof
 We have
 $$F\in K_\bb(I_{X^\vee}) \;\Leftrightarrow\; \xx^{\bb-F} \in I_{X^\vee}
\;\Leftrightarrow\; (\bb-F)^c \in X \;\Leftrightarrow\; F \in
\link(\bb^c,X). \Box $$

\par\noindent
In other words, looking at Betti diagrams (as {\it Macaulay} outputs) 
we have the following picture:
\par\noindent
$\beta_{i,\bb}(I_X):$
$$
\vbox{\offinterlineskip 
\halign{\strut\hfil# \ \vrule\quad&# \ &\hfil # \hfil &\hfil # \hfil &\hfil
# \hfil &\hfil # \hfil
&\hfil # \hfil &\hfil # \hfil &\hfil # \hfil &\hfil # \
\cr
degree&\cr
\noalign {\hrule}
0&1&$h_{-1}(X_{\onedot})$&$h_{-1}(X_{\twodots})$&
$h_{-1}(X_{\threedots})$&$\dots$\cr
1&--&$h_{0}(X_{\twodots})$&$h_{0}(X_{\threedots})$&
$h_{0}(X_{\fourdots})$&$\dots$\cr
2&--&$h_{1}(X_{\threedots})$&$h_{1}(X_{\fourdots})$&
$h_{1}(X_{\fivedots})$&\dots\cr
$\vdots$&\cr
}}$$
where for example $X_{..}$ stands for all full subcomplexes of 
$X$ supported on two vertices, and
\par\noindent
$\beta_{i,\bb^c}(I_{X^\vee}):$
$$
\vbox{\offinterlineskip 
\halign{\strut\hfil$#$ \ \vrule\quad&$#$ \ &\hfil $#$ \hfil &\hfil $#$
\hfil &\hfil $#$ \hfil &\hfil $#$ \hfil
&\hfil $#$ \hfil &\hfil $#$ \hfil &\hfil $#$ \hfil &\hfil $#$ \
\cr
{\rm degree}&\cr
\noalign {\hrule}
0&1&h_{-1}(\link({\onedota}^c))&h_{0}(\link({\twodotsa}^c))
&h_{1}(\link({\threedotsa}^c))&\dots\cr
1&-&h_{-1}(\link({\twodotsa}^c))&h_{0}(\link({\threedotsa}^c))&
h_{1}(\link({\fourdotsa}^c))&\dots\cr
2&-&h_{-1}(\link({\threedotsa}^c))&h_{0}(\link({\fourdotsa}^c))&
h_{1}(\link({\fivedotsa}^c))&\dots\cr
\vdots&\cr
}}
$$
where ${}^c$ stands for complementation, and again the number of dots stands
for the number of vertices in the corresponding faces.

The main observation of this paper is that a simple homological relationship between
restrictions and links has as a consequence the known duality results
involving square-free monomial ideals. We apply it to give a sharper description of
the relationship between the Betti numbers of the dual ideals $I_X$ and $I_{X^\vee}$.

\theorem{thmE}
 For each $\bb \in \set{0,1}^n$ and any vertex $v$ not in
the support of \bb, there is a long exact sequence
 $$\ldots \rightarrow \H_i(X_\bb;k) \rightarrow  \H_i(X_{\bb+v};k) \rightarrow
\H_{i-1}(\link(v,X_{\bb+v});k) \rightarrow \H_{i-1}(X_\bb;k) \rightarrow \ldots.$$

\proof
 This is the long exact homology sequence of the pair $(X_{\bb+v}, X_\bb)$, in
disguise; it is immediate that
 $$\cdots \rightarrow \H_i(X_\bb;k) \rightarrow  \H_i(X_{\bb+v};k) \rightarrow
H_i(X_{\bb+v},X_\bb;k) \rightarrow \H_{i-1}(X_\bb;k) \rightarrow \cdots.$$
 Now, recall that $\star(F,X) = \setdef{G}{F\cup G \in X}$; 
which is the acyclic subcomplex of $X$ generated by all faces
of $X$ which contain $F$. It is also immediate that for all $i$,
 $$H_i(X_{\bb+v}, X_\bb; k) \quad \cong \quad H_i(\star(v,X_{\bb+v}),
\link(v,X_{\bb+v}); k).$$
 Because $\star(v,X_{\bb+v})$ is acyclic, the long exact sequence of the second pair
breaks up into isomorphisms
 $$ H_i(\star(v,X_{\bb+v}), \link(v,X_{\bb+v});k)  \quad \cong \quad 
\H_{i-1}(\link(v,X_{\bb+v});k)$$
 for all $i$. Composing these isomorphisms yields the desired sequence.
 \Box

\ref{thmE} can also be interpreted as the Mayer-Vietoris sequence of the two
subcomplexes $X_\bb$ and $\star(v,X_{\bb+v})$ of $X_{\bb+v}$, whose intersection
is $\link(v,X_{\bb+v})$.

We shall exploit the exactness of this sequence at $\H_i(X_\bb)$. It is easy to observe
this exactness at the level of cycles: Let $\alpha$ be an $i$-cycle supported on
$X_\bb$, representing a homology class in $\H_i(X_\bb;k)$. If $\alpha$ maps to zero
in $\H_i(X_{\bb+v};k)$, then there exist an $(i+1)$-cycle $\beta$ supported on
$X_{\bb+v}$, whose boundary $\partial\beta = \alpha$. Express $\beta$ as a sum
$\beta_1+\beta_2$, where $\beta_1$ is supported on $X_\bb$ and every face of
$\beta_2$ contains the vertex $v$. Define $\alpha'=\partial\beta_2 = \alpha -
\partial\beta_1$. The cycle $\alpha'$ is supported on  $\link(v,X_{\bb+v})$, and
represents the same homology class as $\alpha$ in $\H_i(X_\bb;k)$.

\corollary{corF}
 The Betti numbers of $I_X$ and of $I_{X^\vee}$ satisfy the inequality
 $$\beta_{i,\bb} \quad\le\quad \sum_{\bb\, \preceq \,\cc\, \preceq \,[n]}
\; \beta^\vee_{\abs\bb-i-1,\cc}$$
 for each $0\le i \le n-1$ and each $\bb \in \set{0,1}^n$.

\proof
 The exactness at $\H_i(X_\bb;k)$ of the sequence of \ref{thmE} yields the inequality
 $$ \dim \H_i(X_\bb;k) \quad \le \quad \dim \H_i(\link(v,X_{\bb+v});k) \; + \;
\dim \H_i(X_{\bb+v};k).$$
 Note that for any face $F$ disjoint from \bb, and any vertex $v$ not in the support
of $\bb+F$,
 $$\link(v,\link(F,X)_{\bb+F+v}) \quad = \quad \link(F+v,X_{\bb+F+v}).$$ Applying
\ref{thmE} to ${\link(F,X)}_{\bb+F}$ in place of $X_\bb$ yields the exact sequence
 $$\H_i(\link(F+v,X_{\bb+F+v});k) \longrightarrow \H_i(\link(F,X_{\bb+F});k)
\longrightarrow  \H_i(\link(F,X_{\bb+F+v});k).$$
Combining the resulting inequalities while iteratively adding vertices  yields
 $$ \dim \H_i(X_\bb;k) \quad \le \quad \sum_{F \,\cap\, \bb\, =\, \emptyset} \dim
\H_i(\link(F,X);k).$$
 By \refs{thmC} and \refn{thmD} these dimensions can be interpreted as Betti
numbers of $I_X$ and $I_{X^\vee}$, respectively.
 \Box

In particular, summing up and collecting all terms of the same total
degree we obtain:

\corollary{corG} The single graded Betti numbers of $I_X$ and of $I_{X^\vee}$ 
satisfy the inequality
$$\beta_{i,m} \quad\le\quad \sum_{k=0}^{n-m}\;{{m+k}\choose k}
\beta^\vee_{m-i-1,m+k},$$
for each $0\le i \le n-1$ and each $m\ge i+1$.

The following consequence of \ref{thmE} and \ref{corF} extends Terai's
characterization of dual Stanley-Reisner ideals. Define a Betti number
$\beta_{i,\bb}$ to be {\it $i$-extremal\/} if $\beta_{i,\cc}=0$ for all $\cc\succ\bb$,
that is  all multigraded entries below $\bb$ on the $i$-th column vanish
in the Betti diagram as a {\it Macaulay} output [Mac].
Define $\beta_{i,\bb}$ to be {\it extremal\/} if $\beta_{j,\cc}=0$ for all
$j \ge i$, and $\cc\succ\bb$ so $\abs\cc-\abs\bb \ge j-i$. In other words,
$\beta_{i,\bb}$ corresponds to the ``top left corner'' of a box of zeroes
in the multigraded Betti diagram, thus our definition agrees with the single graded
one we've introduced in \ref{gins}.  Note that we have not
assumed this time that $\beta_{i,\bb} \ne 0$.

\theorem{thmG}
 If $\beta^\vee_{i,\bb}$ is $i$-extremal, then the inclusion $\link(\bb^c,X) \subseteq
X_b$ induces an exact sequence
 $$\H_{i}(X_\bb,\link(\bb^c,X); k)\longrightarrow
\H_{i-1}(\link(\bb^c,X); k) \longrightarrow \H_{i-1}(X_\bb; k) \longrightarrow 0$$
showing that 
$$\beta^\vee_{i,\bb} \ge \beta_{\abs\bb-i-1,\bb}.$$ If
$\beta^\vee_{i,\bb}$ is extremal, then the above surjection is in fact an isomorphism, 
showing that
$\beta^\vee_{i,\bb} = \beta_{\abs\bb-i-1,\bb}$.

\proof
The condition that $\beta^\vee_{i,\bb}$ is $i$-extremal, 
means that the right hand sum in \ref{corF}, applied for $\abs\bb-i-1$ instead of $i$,
has exactly one summand, which gives the first part of the theorem. If
moreover $\beta^\vee_{i,\bb}$ is extremal, then $\beta_{\abs\bb-i-1,\bb}$ is
in fact $(\abs\bb-i-1)$-extremal so the second claim follows from the first part
applied for $X^\vee$ instead of $X$. 
\Box

Looking at a Betti diagram as output by {\it Macaulay}, this result asserts in
particular that any lower right corner flips via duality. Thus we can speak of
``$d$-regularity in homological dimensions $\ge i$'' and interpret it as a
 statement generalizing Terai's theorem [Te97] 
(compare [FT97, Corollary 3.2]):

\corollary{terai}
 The regularity of $I_X$  equals  the projective dimension
of $S/I_{X^\vee}$, the Stanley-Reisner ring of its Alexander dual.

\proof The regularity of $I_X$ is computed by the largest 
$\abs\bb -i$ such that $\beta_{i,\bb}\ne 0$, while the projective
dimension of $S/I_{X^\vee}$, is $1+\max{(j)}$ such that $\beta^\vee_{j,\cc}
\ne 0$ for some $\cc\in \set{0,1}^n$. Thus the claim follows from \ref{thmG},
because of the equality of the corresponding pairs of extremal Betti numbers.
\Box

Moreover, \ref{thmG} provides also easy proofs of classical
criteria due to Reisner [Re76], and Stanley [Sta77] respectively:

\theorem{CM}
The following conditions are equivalent:
\item{a)} $S/I_X$ is a Cohen-Macaulay ring;
\item{b)} $\H_{i}(\link(F,X); k)=0$,  for all $F\in X$ and $i<\dim(\link(F,X))$.

\proof
By \ref{thmG} or \ref{terai}, if $S/I_X$ is Cohen-Macaulay, 
then $I_{X^\vee}$ is generated in
degree $n-\dim(X)-1$ and  has a linear resolution.
In other words $\beta^\vee_{i,\bb}=0$,
for all $\bb$ with $\abs\bb > n-\dim(X)-1+i$. By \ref{thmD},
this means that for all $F=\bb^c\in X$, 
$\dim \H_{i}(\link(F,X); k)=0$ for $i< \dim(X)+\abs\bb -n=
\dim(X)-\abs F= \dim(\link(F,X))$. To prove the implication
$b)\Rightarrow a)$ it is enough to show that $X$ is pure,
and then the above argument reverses.
Since $\link(G,\link(F,X))=\link(F\cup G, X)$, whenever $F\cup G\in X$
and $F\cap G=\emptyset$, we observe that the same cohomological
vanishing holds for all proper links of $X$, hence by induction
we may assume that they are pure. Now if $\dim(X)\ge 1$, then
$b)$ also gives $\H_{0}(X;k)=\H_{0}(\link(\emptyset,X); k)=0$,
so $X$ is connected and this together with the purity of the links
shows the purity of $X$.
\Box

Since $S/I_X=S/I_{\core(X)}[X_i\mid v_i\in X\setminus \core(X)]$
(see for instance [BH93, p.232]), we have
that $S/I_X$ is Gorenstein iff $S/I_{\core(X)}$ is Gorenstein, thus it
is enough to show the following

\theorem{Gor}
If $X=\core(X)$ (that is $X$ is not a cone), 
then the following are equivalent:
\item{a)} $S/I_X$ is a Gorenstein ring (over $k$);
\item{b)} For all $F\in X$, $\H_{i}(\link(F,X); k)=
\cases{$k$ & if $i=\dim(\link(F,X))$,\cr 0 & otherwise}$;

\proof To prove the implication $a)\Rightarrow b)$, 
we argue by induction on $\abs F$:
for any vertex $v\in F$, one has 
$S/I_{\link(v,X)}=S/(I_X:(x_v))$,
on the other hand $S/I_X$ Gorenstein implies that 
$S/(I_X:(x_v))=x_vS/I_X$
is also Gorenstein, whereas 
$\link(G,\link(v,X))=\link(v+G,X)$, for all $\{v\}\cup G\in X$
with $v\not\in G$.

 If condition $b)$ holds, then $S/I_X$ is a Cohen-Macaulay ring
by \ref{CM}, and so $X$ is pure. Moreover $X$ is
a pseudomanifold, that is every $(\dim(X)-1)$-face
of $X$ lies in exactly two facets, and $X$ is orientable,  
that is $\H_{\dim(X)}(X; k)=k$ since $X=\link(\emptyset,X)$.
In fact the same holds for every proper link of $X$. We
observe next that $X_{v^c}$ is also Cohen-Macaulay of the same
dimension, for any vertex $v$ of $X$. By \ref{CM}, all we have to 
check is that  $\H_{i}(\link(F,X_{v^c}); k)=0$,  for all $F\in X_{v^c}$ and 
$i<\dim(\link(F,X_{v^c}))$. If $v$ is not a vertex of
$\link(F,X)$, or $i<\dim(\link(F,X_{v^c}))-1$ this is immediate
from condition $b)$. If $v$ is a vertex of $\link(F,X)$ and 
$i=\dim(\link(F,X_{v^c}))-1$, this follows from the long exact
sequence in \ref{thmE}
$$\rightarrow \H_{i+1}(\link(F,X);k) \rightarrow  
\H_i(\link(F\cup\{v\},X);k) \rightarrow
\H_i(\link(F,X_{v^c});k) \rightarrow 
\H_i(\link(F,X);k)$$
where the leftmost arrow is nonzero being
induced by the restriction of an orientation class.
To prove that $S/I_X$ is Gorenstein, 
it is enough to show that
the canonical module of $S/I_X$ is invertible, or equivalently that
all generators of the canonical module lie in a single degree and
$S/I_X$ has a unique extremal Betti number whose value is one.
The first condition follows now from the fact 
that $X_{v^c}$ is also Cohen-Macaulay of the same
dimension, for any vertex $v$ of $X$, and thus 
$\beta_{n-\dim(X)-2,\bb}=0$, for all
$\bb\ne(1,1,\ldots,1)$. For the second condition observe that,
by \ref{terai}, $I_{X^\vee}$ has a linear resolution
and thus by \ref{thmG} the unique extremal Betti number
of $I_{X^\vee}$, which is one by our hypothesis,
coincides with the corresponding extremal Betti number
of $I_X$.\Box

\remark{dCM} Recall that a simplicial complex $X$ is called {\it doubly
Cohen-Macaulay} if $X$ is a Cohen-Macaulay complex 
and $X_{v^c}$ is also Cohen-Macaulay of the same
dimension, for any vertex $v$ of $X$. The formula in \ref{thmD} and
the proof of \ref{Gor} show that if $X$ is doubly
Cohen-Macaulay, then $I_{X^\vee}$ has a linear resolution and
$\beta^\vee_{\dim(X)+1,(1,1,\ldots,1)}$ is the unique extremal 
Betti number.

\section{examples} Examples

We end with three examples illustrating the above described behavior of the 
extremal multigraded Betti numbers:

\example{cycle}
Let $X$ be  a length five cycle, that is 
$I_X={(x_ix_{i+2})}_{i\in\ZZ_5}\subset 
k[x_0,\ldots,x_4]$. Then $X^\vee$ is the triangulation of a 
M\"obius band shown in
Figure~1,  and $I_{X^\vee}={(x_ix_{i+1}x_{i+2})}_{i\in\ZZ_5}$.

$$
\vbox{\offinterlineskip 
\halign{\strut\hfil# \ \vrule\quad&# \ &\hfil # \ &\hfil # \ &\hfil # \ &\hfil # \ 
&\hfil # \   
\cr
degree&1&5&5&1\cr
\noalign{\hrule}
0&1&--&--&--\cr
1&--&\ 5&\ 5&--\cr
2&--&--&--&\ 1\cr
\noalign{\bigskip}
\omit&\multispan{3}{$\beta_{i,j}(I_X)$}\cr
\noalign{\smallskip}
}}\hskip 2.2truecm
\vbox{\offinterlineskip 
\halign{\strut\hfil# \ \vrule\quad&# \ &\hfil # \ &\hfil # \ &\hfil # \ &\hfil # \ 
&\hfil # \ 
\cr
degree&1&5&5&1\cr
\noalign {\hrule}
0&1&--&--&--\cr
1&--&--&--&--\cr
2&--&\ 5&\ 5&\ 1\cr
\noalign{\bigskip}
\omit&\multispan{3}{$\beta_{i,j}(I_{X^\vee})$}\cr
\noalign{\smallskip}
}}
$$

\draw{75}{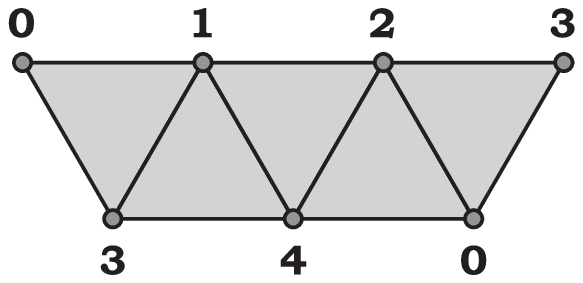}{}

\example{torus}
It is easily seen that a  triangulation of the torus $T_1$ has at least
7 vertices, and in case the triangulation has exactly 7 vertices, 
that the graph of its 1-skeleton is necessarily $K_7$,
the complete graph on seven vertices.
Such a triangulation $X$ (first constructed in 
1949 by Cs\'asz\'ar) is shown in Figure~2; it is unique
up to isomorphism and has an automorphism group of order 42.
The dual graph of its 1-skeleton divides the torus in 
the well known 7-colourable map (see [Wh] for more details).
Thus up to a permutation, $I_X={(x_ix_{i+1}x_{i+2},\, x_ix_{i+1}x_{i+4},\,  
x_ix_{i+2}x_{i+4})}_{i\in\ZZ_7}$.  Then 
$I_{X^\vee}={(x_ix_{i+1}x_{i+2}x_{i+4},\, x_ix_{i+1}x_{i+2}x_{i+5})}_{i\in\ZZ_7}$.

$$
\vbox{\offinterlineskip 
\halign{\strut\hfil# \ \vrule\quad&# \ &\hfil # \ &\hfil # \ &\hfil # \ &\hfil # \ 
&\hfil # \ &\hfil # \ &\hfil # \ &\hfil # \ &\hfil # \  
\cr
degree&1&21&49&42&15&2\cr
\noalign {\hrule}
0&1&--&--&--&--&--\cr
1&--&--&--&--&--&--\cr
2&--&21&49&42&14&\ 2\cr
3&--&--&--&--&\hfill 1&--\cr
\noalign{\bigskip}
\omit&\multispan{5}{$\beta_{i,j}(I_X)$}\cr
\noalign{\smallskip}
}}
\hskip 2.2truecm
\vbox{\offinterlineskip 
\halign{\strut\hfil# \ \vrule\quad&# \ &\hfil # \ &\hfil # \ &\hfil # \ &\hfil # \ 
&\hfil # \ &\hfil # \ &\hfil # \ &\hfil # \ &\hfil # \ 
\cr
degree&1&14&21&9&1\cr
\noalign {\hrule}
0&1&--&--&--&--\cr
1&--&--&--&--&--\cr
2&--&--&--&--&--\cr
3&--&14&21&\ 7&\ 1\cr
4&--&--&--&\ 2&--\cr
\noalign{\bigskip}
\omit&\multispan{5}{$\beta_{i,j}(I_{X^\vee})$}\cr
\noalign{\smallskip}
}}
$$

\draw{75}{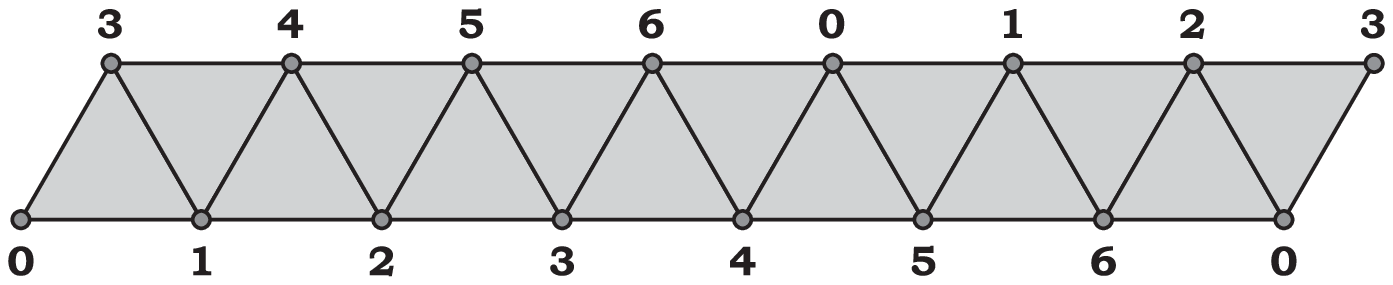}{}

\example{brunsification} In fact, one can  construct examples 
of homogeneous modules with prescribed extremal graded Betti numbers,
for example, by considering appropriate
direct sums where each direct summand features exactly one 
extremal Betti number. Moreover, a classical result
of Bruns [Br76] (see also [EG85, Corollary 3.13, p.56]) 
implies that all such possible extremal ``shapes'' 
and values of extremal Betti numbers in resolutions of modules may
be realized also in minimal free resolutions of homogeneous ideals 
(generated by 3 elements). By passing to the generic initial ideal
and then polarizing we may also construct examples of squarefree
monomial ideals with the desired extremal Betti numbers.

\example{multigraded} Extremal multigraded numbers need not to
be also extremal in the total degree sense. For example, if $X$ is the simplicial
complex shown in Figure~3, then $I_X=(x_0x_2, x_0x_3, x_0x_4, x_1x_4)$,
and $I_{X^\vee}=(x_0x_4, x_0x_1, x_2x_3x_4)$, 

\draw{85}{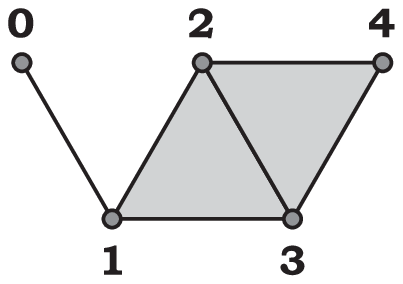}{}

\noindent
while the corresponding Betti diagrams are

$$
\vbox{\offinterlineskip 
\halign{\strut\hfil# \ \vrule\quad&# \ &\hfil # \ &\hfil # \ &\hfil # \ &\hfil # \ 
&\hfil # \   
\cr
degree&1&4&4&1\cr
\noalign{\hrule}
0&1&--&--&--\cr
1&--&\ 4&\ 4&\ 1\cr
\noalign{\bigskip}
\omit&\multispan{3}{$\beta_{i,j}(I_X)$}\cr
\noalign{\smallskip}
}}\hskip 2.2truecm
\vbox{\offinterlineskip 
\halign{\strut\hfil# \ \vrule\quad&# \ &\hfil # \ &\hfil # \ &\hfil # \ &\hfil # \  
\cr
degree&1&3&2\cr
\noalign {\hrule}
0&1&--&--\cr
1&--&\ 2&\ 1\cr
2&--&\ 1&\ 1\cr
\noalign{\bigskip}
\omit&\multispan{3}{$\beta_{i,j}(I_{X^\vee})$}\cr
\noalign{\smallskip}
}}
$$

Both second order syzygies of 
$I_{X^\vee}$ are extremal (in the multigraded sense),
but $\beta_{1,\{0,1,4\}}=\beta^\vee_{1,\{0,1,4\}}$, which is 
also extremal, is not
extremal in the single graded sense.

\references

\item{[Bay96]} Bayer, D., Monomial ideals and duality, 
Lecture notes, Berkeley 1995-96, available by anonymous ftp from \hfill\break
{\tt ftp://math.columbia.edu/pub/bayer/monomials\_duality/}.

\item{[BaSt87]} Bayer, D., Stillman, M., 
A criterion for detecting $m$-regularity. 
{\it Invent. Math.} {\bf 87}, (1987), no. 1, 1--11.

\item{[Mac]} Bayer, D., Stillman, M.,
Macaulay: A system for computation in
        algebraic geometry and commutative algebra.
Source and object code available for Unix and Macintosh
computers. Contact the authors, or download 
from {\tt ftp://math.harvard.edu/}.

\item{[Br76]} Bruns, W., ``Jede'' endliche freie Aufl\"osung ist freie Aufl\"osung eines von drei Elementen erzeugten Ideals. 
{\it J. Algebra} {\bf 39}, (1976), no. 2, 429--439.

\item{[BH93]} Bruns, W., Herzog, J., {\it Cohen-Macaulay Rings},  Cambridge
Studies in advanced mathematics, {\bf 39}, Cambridge University Press 1993.

\item{[Ea95]} Eagon, J., Minimal resolutions of ideals associated to homology
triangulated manifolds, preprint 1995.

\item{[ER96]} Eagon, J., Reiner, V., Resolutions of Stanley-Reisner rings
and Alexander duality, preprint 1996.

\item{[Ei95]} Eisenbud, D.,
{\it Commutative Algebra with a View Toward Algebraic Geometry}, Springer, 
1995.

\item{[EG85]} Evans E.G., Griffith, P., {\it Syzygies}, 
London Mathematical Society LNS, {\bf 106}, 
Cambridge University Press, 1985.

\item{[FT97]} Fr\"uhbis, A., Terai, N., Bounds for the regularity of
monomial ideals, preprint 1997

\item{[Ho77]} Hochster, M.,Cohen-Macaulay rings, combinatorics and simplicial
complexes, in Ring theory II,  McDonald B.R., Morris, R. A. (eds),
{\it Lecture Notes in Pure and Appl. Math.}, {\bf 26}, M. Dekker 1977.

\item{[Mun84]} Munkres, J. R., {\it Elements of algebraic topology},
Benjamin/Cummings, Menlo Park CA, 1984.

\item{[Rei76]} Reisner, G.A., Cohen-Macaulay quotients of polynomial rings,
{\it Adv. in Math.}, {\bf 21}, (1976) 30--49.

\item{[Sta77]} Stanley, R., Cohen-Macaulay complexes, in {\it Higher
Combinatorics}, (M. Aigner, ed.), Reidel, Dordrecht and Boston, 1977.

\item{[Sta96]} Stanley, R., Combinatorics and Commutative Algebra, Second
edition, Progress in Math. {\bf 41}, Birkh\"auser, 1996.

\item{[Te97]} Terai, N., Generalization of Eagon-Reiner theorem
and $h$-vectors of graded rings, preprint 1997.

\item{[Wh]} White, A., Graphs, Groups and Surfaces,
{\it North-Holland Mathematics Studies}, {\bf 8}, North Holland 1984.

\bigskip\bigskip
\vbox{\noindent Author Addresses:
\smallskip
\noindent{Dave Bayer}\par
\noindent{Department of Mathematics, Barnard College, Columbia University, 
2990 Broadway, New York, NY 10027}\par
\noindent{bayer@math.columbia.edu}
\medskip
\noindent{Hara Charalambous}\par
\noindent{Department of Mathematics,  University at Albany, SUNY,
Albany, NY 12222}\par
\noindent{hara@math.albany.edu}
\medskip
\noindent{Sorin Popescu}\par
\noindent{Department of Mathematics, Columbia University,
2990 Broadway, New  York, NY 10027}\par
\noindent{psorin@math.columbia.edu}\par
}

\bye